\documentclass[12pt]{amsart}
\usepackage{epsfig,amssymb,epic,eepic,color}
\usepackage[all]{xy}
 \usepackage[T1]{fontenc}        
\newtheorem{thm}{Theorem}[section]

\newtheorem{prop}[thm]{Proposition}

\newtheorem{lemme}[thm]{Lemma}
\newtheorem{cor}[thm]{Corollary}

\newtheorem{remarque}[thm]{Remark}

\newtheorem{rien}[thm]{}

\newcommand{\be}{\begin{enumerate}}
\newcommand{\ee}{\end{enumerate}}
\newcommand{\bi}{\begin{itemize}}
\newcommand{\ei}{\end{itemize}}

\def\R{\mathbb{R}}

\def\Z{\mathbb{Z}}

\def\T{\mathbb{T}}

\def\ga{\gamma}

\def\al{\alpha}
\def\be{\beta}

\def\De{\Delta}

\def\si{\sigma}
\def\Si{\Sigma}

\def\nd{\noindent}
\def\bull{\hfill$\Box$}
\def\proof{\nd {\bf Proof.\ }}

\begin{document}
\vskip 1cm
\title{Morse 2-jet space and $h$-principle}

%\author{Alain CHENCINER}

%\author{Fran\c cois LAUDENBACH}

%\begin{center}{\sc }
%\end{center}
%\vskip 1cm
%\begin{center}
%{\sc }
%\end{center}

\author{Alain Chenciner}
\address{Universit\'e Paris 7 et IMCCE (UMR 8025 du CNRS)}
\email{chenciner@imcce.fr}

\author{Fran\c cois Laudenbach}
\address{Universit\'e de Nantes, UMR 6629 du CNRS, 44322 Nantes, France}
\email{francois.laudenbach@univ-nantes.fr}

\date{This version: February, 2009.}

\keywords{Morse singularity, 2-jet, $h$-principle}

\subjclass[2000] {58A20, 58E05, 57R45}

\thanks{}

\begin{abstract}
A section in the $2$-jet space of Morse functions is not always homotopic to a holonomic section. We give a necessary  condition for being the case and we discuss
the sufficiency.
\end{abstract}

\maketitle

\thispagestyle{empty}
\vskip 1cm
\section{Introduction}
\medskip

Given a submanifold $\Sigma$ in an $r$-jet space (of smooth sections of a bundle over a manifold $M$), it is natural to look at the associated {\it differential relation} 
$\mathcal R(\Si)$ formed by the $(r+1)$-jets  {\it transverse} to $\Si$. 
%This makes sense 
For $j^rf$ being  transverse to $\Si$ at $x\in M$ is detected by $j^{r+1}_xf$. This is an open differential relation in the corresponding $(r+1)$-jet space.
One can ask whether the Gromov {\it  $h$-principle} holds: is any section with value in 
$\mathcal R(\Si)$ homotopic to a {\it holonomic} section of $\mathcal R(\Si)$? (We recall that a holonomic section of a $(r+1)$-jet space is  a section of the form $j^{r+1}f$.)

According to M. Gromov, the answer is yes when $M$ is an open manifold and 
$\Sigma$ is {\it natural}, that is, invariant by a lift  of $Diff(M)$ to the considered jet space (see \cite{gromov} p. 79, \cite{el-mi} ch. 7). 

%When looking at jet space of functions, 
The answer is also yes
when the codimension of $\Si$ is %bigger  
higher than the dimension $n$ of $M$; % in the case of jet space of functions 
 this case follows easily from Thom's transversality theorem in jet spaces (see \cite{bordeaux}). In the case of jet space of functions and when $\Si$ is natural and 
 ${\rm codim}\Si\geq n+1$, 
 it also can be seen as a baby case of a theorem of Vassiliev \cite{vassiliev}. 
 %%Let $D$ (like {\it discriminant}) denote the complement of $\mathcal R(\Si)$ in the $(r+1)$-jet space of functions $J^{r+1}(M)=J^{r+1}(M,\R)$. As a generic section of $J^r(M)$ avoids $\Si$, $\mathcal R(\Si)$ consists of $(r+1)$-jets whose $r$-jet lies in the complement of $\Si$ in $J^r(M)$. Hence ${\rm codim} D\geq n+1$. In that case Vassiliev's theorem states  that the  $h$-principle holds for the complement of $D$, even with parameters
 
% Here we consider jets of functions, $r=1$
In this note we are interested in a codimension $n$ case when $M$ is a compact $n$-dimensional manifold. % the case of $2$-jets of Morse functions. 
Let $J^r(M)$ denote the space of $r$-jets of real functions; when the boundary of $M$ is not empty, it is meant that we speak of jets of functions which are locally constant on the boundary. We take 
$\Sigma \subset J^1(M)$
the set of critical $1$-jets. Then $\mathcal R(\Si)\subset J^2(M)$ is the open set 
of $2$-jets of Morse functions. We shall analyze the obstructions preventing the  $h$-principle to hold with this differential relation.\\

\section{Index cocycles}

\medskip
It is more convenient to work with the {\it reduced} jet spaces $\widetilde J^r(M)$, quotient of $J^r(M)$ by $\R$ which acts by translating the value of the jet. It is a vector bundle whose linear structure is induced by that of $C^\infty(M)$. For instance,
 $\widetilde J^1(M)$ is isomorphic to the cotangent space $T^*M$. Let
$\mathcal M$ denote the reduced $2$-jets of Morse functions, that is the $2$-jets which are transverse to the zero section $0_M$ of $T^*M$ (in the sequel, {\it jet } will mean {\it reduced jet}). Let $\pi: \widetilde J^2(M)\to \widetilde J^1(M)$ be the projection 
and $\pi_0:\widetilde J^2_0(M)\to 0_M$ be its restriction over the zero section 
of the cotangent space. Since it is formed of critical $2$-jets, it is a vector bundle whose fiber is the space of quadratic forms, $S^2(T^*_xM)$, $x\in M$.
 Let $\mathcal M_0:= \mathcal M\cap \widetilde J^2_0(M)$; it is a bundle  over $0_M$ whose fiber consists of 
 non-degenerate quadratic forms.  Its complement in $\widetilde J^2_0(M)$ is denoted 
 $\mathcal D$ (like discriminant); it is formed of $2$-jets which are not transverse to 
 $0_M$. When $M$ is connected,  $\mathcal M_0$ has a connected component 
 $\mathcal M_0^i$ for each index $i\in\{1, ..., n\}$ of quadratic forms. \\
 
 \begin{rien} Tranverse orientation. {\rm Each $\mathcal M_0^i$ is a proper submanifold of codimension $n$ in $\mathcal M$. Moreover the differential $d\pi$ gives rise to 
 an isomorphism of normal fiber bundles} 
 $$\nu(\mathcal M_0^i, \mathcal M)\cong \pi^*(\nu(0_M,T^*M))\vert \mathcal M_0^i\,.
 $$
 \end{rien}
 Of course, $\nu(0_M,T^*M)$ is canonically isomorphic to the cotangent bundle $\tau^*M$, whose total space is $T^*M$. When $M$ is oriented,  so is the bundle 
 $\tau^*M$. When $M$ is not orientable, one has a local system of orientations of 
 $\tau^*M$. Pulling it back by $\pi$ yields a local system of orientations 
 of $\nu(\mathcal M_0^i, \mathcal M)$ (that is, co-orientations of $\mathcal M_0^i$). Let us denote  $\mathcal M_0^{even}$ (resp.
 $\mathcal M_0^{odd}$) the union of the $\mathcal M_0^{i}$'s for $i$ even (resp. odd).
 We endow $\mathcal M_0^{even}$ with the above local system of co-orientations. For reasons
 which clearly appear below, it is more natural to equip $\mathcal M_0^{odd}$ with the opposite system of co-orientations.\\
 
 \begin{lemme}\label{positive} Let $s=j^2f$ be a %(local) 
 holonomic section of $\mathcal M$ meeting 
 $\mathcal M_0$ transversally. Then 
 each intersection point of $s(M)$ with $\mathcal M_0$ is positive. The same statement holds when $s$ is  a local holonomic  section only.
 \end{lemme}
 
 \proof  Let $a$ be such an intersection point in $s(M)\cap \mathcal M_0^i$; so $i$ is the index of the corresponding critical $2$-jet. We can calculate in local coordinates 
 $(x,y',y'')$, where $x=(x_1,...,x_n)$ are local coordinates of $M$,% with $x(a)=0$,
 $y'=(y'_1,...,y'_n)$ (resp. $y''=(y''_{jk})_{1\leq j\leq k\leq n}$) are the associated coordinates of $T^*_xM$ (resp. $S^2T^*_xM$). Since $f$ is holonomic, we have 
 $y''_{jk}(a)= \frac{\partial y'_j}{\partial x_k}(a)$.  Finally, the sign of 
 $\det y''(a)$ (positive if $i$ is even and negative if not) gives the sign of the Jacobian determinant at $a$ of the map $x\mapsto y'(x)$, that is the sign of the intersection point 
 when $\mathcal M^i_0$ is co-oriented by the canonical orientation of the $y'$-space.
 As we have reversed this co-orientation when $i$ is odd, the intersection point is positive whatever the index is.\bull\\

 \begin{prop}\label{ineq}
 
 \nd {\rm 1)}
 Each  $\mathcal M_0^i$ defines a degree $n$ cocycle of $\mathcal M$ with coefficients in the local system
 $\Z^{or}$ of integers twisted by the orientation of $M$. Let $\mu_i$ be its cohomology
 class in $H^n(\mathcal M, \Z^{or})$; in particular, if $s:M\to\mathcal M$ is a section,
 $<\mu_i,[s]>$ is an integer.
 
\nd {\rm  2)} When $s$ is homotopic to a holonomic  section $j^2f$, then  $<\mu_i,[s]>$ is positive and equals the number $c_i(f)$ of critical points of the Morse function $f$. In particular the total number $\vert Z\vert$ of zeroes of the section $\pi\circ s$ (which, by construction, is transverse to the 0-section) satisfies:
$$\vert Z\vert\geq \sum_{i=0}^n c_i(f)\,.$$
\end{prop}
 
 \proof 1) Let $\sigma$ be a singular  $n$-cycle with twisted coefficients of 
 $\mathcal M$. It can be $C^0$-approximated by $\si'$, an $n$-cycle which is transverse to $\mathcal M_0^i$.
 As $\mathcal M_0^i$ is a proper submanifold, there are finitely many 
 intersections points  in $\si'\cap \mathcal M_0^i$, each one having a sign with respect to the local system of coefficients. The algebraic sum of these signs defines an integer
$c(\sigma')$. One easily checks that $c(\sigma')=0$ if $\sigma'$ is a boundary. As a consequence, if $\si'_0$ and $\si'_1$ are two approximations of $\si$, as 
$\si'_1-\si'_0$ is a boundary, we have $c(\si'_1)-c(\si'_0)=0$
 which allows us to uniquely define 
$c(\si)$ as the value of an $n$-cocycle  on $\si$.  Typically, the image of a section carries an $n$-cycle with twisted coefficients and this algebraic counting applies.

\nd 2) Since  $c$ defined in 1) is a cocycle, it takes the same value on $s$ and on 
$j^2f$. According to  lemma \ref{positive}, it counts +1 for each intersection point
in $j^2f\cap \mathcal M^i_0$, that is, for each index $i$ critical point of $f$.\bull\\

\begin{cor}\label{euler}If $s$ is a section of $\mathcal M$ which is homotopic to a holonomic section, the integers $m_i:=<\mu_i, [s]>$ fulfill the Morse inequalities
$$\begin{array}{rcl}m_0 &\geq &\beta_0(F)\\
m_1-m_0 &\geq &\beta_1(F)-\beta_0(F)\\
... & & ...\\
m_0 -m_1+...+ (-1)^n m_n&=&\beta_0(F)-\beta_1(F) ...+(-1)^n\beta_n(F)= :\chi(M)
\end{array}
$$where $F$ is a field of coefficients,  $\beta_i(F)=\dim_FH_i(M,F^{or})$ is the $i$-th
Betti number with coeffcients in $F^{or}$ ($F$ twisted by the orientation)
and $\chi(M)$ 
is the Euler characteristic (independent of the field $F$).\\
\end{cor}

\begin{cor} The  $h$-principle does not hold true for the sections of $\mathcal M$.
\end{cor} 

\proof   It is sufficient to construct a section $s$ of $\mathcal M$ which violates
the Morse inequalities, for example a section which does not intersect $\mathcal M_0^0$. Leaving the case of the circle as an exercise, we may assume 
$n>1$.
One starts with a section $s_1$ of $T^*M$ tranverse to $O_M$. Each zero 
of $s_1$
has a sign (if the local orientation of $M$
 is changed, so are both local orientations of $s_1$ and $0_M$ the sign of the zero in unchanged). For each zero $a$, one can construct a homotopy fixing $a$, with arbitrary small support, which makes $s_1$ linear in a small neighborhood of $a$. As $GL(n,\R)$ has exactly two connected components, one can even suppose that after the homotopy, $s_1$ is near $a$ the derivative of a non degenerate quadratic fonction whose index can be chosen arbitrarily provided it is even (resp. odd) if $a$ is a positive (resp. a negative) zero. Finally, one can achieve by homotopy that near each zero $a$, one has $s_1=df$ with $a$  a
 non-degenerate critical point of $f$ of  index 2 (resp. 1)  if $a$ is a positive (resp. negative) zero. 
 
 Near the zeroes $s_1$ has a canonical lift to $\mathcal M$ by $s_2= j^2f$. Away from the zeroes, the lift $s_2$ extends as a lift of $s_1$ since  the fibers of $\pi $
are contractible over $T^*M\setminus 0_M$. By construction, we have
$<\mu_0, [s_2]>=0$, violating the first Morse inequality.\bull\\

\begin{remarque} Denote $\mu_{even}= \mu_0 +\mu_2+...$ and  $\mu_{odd}= \mu_1 
+...$\,. The following statement holds true: $\mu_{even}=\mu_{odd}$ if and only if 
the Euler characteristic vanishes.
\end{remarque}

\proof Assume first $\mu_{even}=\mu_{odd}$. Proposition \ref{ineq} yields for any holonomic section in $\mathcal M$:
$ m_{even}=m_{odd}$, that is $\chi(M)=0$. Conversely, if $\chi(M)=0$, there exists a 
non-vanishing $1$-form on $M$ and hence, by lifting it to $\widetilde J^2(M)$, a section $v_0$ in $\mathcal M$ avoiding
 $\mathcal M_0$. We form 
 $$W=\{z\in \widetilde J^2(M) \mid z=z_0+tv_0, z_0\in \mathcal M_0,t\geq 0
 \ {\rm or}\ z_0\in \mathcal D, t>0\}.$$
 It is a proper submanifold in $\mathcal M$ whose boundary (with orientation twisted coefficients) is $\mathcal M_0^{even}-\mathcal M_0^{odd}$. Therefore, every cycle $c$ satisfies $<\mu_{even},c>=<\mu_{odd},c>$, which implies  the wanted equality. \bull%$\mu_{even}=\mu_{odd}$\,.\bull
 \goodbreak

\section{ Are Morse inequalities  sufficient?}

\medskip
This question is closely related to the problem of minimizing the number of critical points of a Morse function. This problem  was solved by S. Smale in dimension higher than 5 for simply connected manifolds, as a %the
consequence of the methods he developped
for proving his famous $h$-cobordism theorem (see \cite{smale} or chapter 2 in \cite{franks}). Under the same topological assumptions %as Smale 
we can answer our question positively. But there are other cases, discussed later, where the answer is negative.

\begin{prop} Two sections $s,t$ of $\mathcal M\subset\widetilde J^2(M)$ are homotopic as sections of 
 $\mathcal M$ if and only if their algebraic intersection numbers  $m_i$ with $\mathcal M_0^i$ are the same.\label{homot}
\end{prop}

\proof According to proposition \ref{ineq} 1), the condition is necessary. Let us prove that it is sufficient. Leaving the 1-dimensional case to the reader, we assume $\dim M\geq 2$.
Denote $s^1=\pi\circ s$. Each zero of $s^1$ is given an index due to its lifting by $s$ to a point of some $\mathcal M_0^i$. 
For each index $i$ choose $|m_i|$ zeroes of $s^1$, 
$a_i^1,\ldots, a_i^{|m_i|}$, among its zeroes of 
index $i$; % (corresponding to intersection points of $s$ with $\mathcal M_0^i$); 
when $m_i>0$ (resp. $m_i<0$),  we choose the $a_i^j$ so that the corresponding intersection points of
$s(M)$ with ${\mathcal M}_0^i$ are positive (resp. negative).  
When $m_i=0$, no points are selected. In the same way, $|m_i|$ zeros $b_i^1,\ldots, b_i^{|m_i|}$ of $t^1$ are chosen.

The intersection signs being the same,  one can find a homotopy of $t$ in $\mathcal M$, which brings the $b_i^j$ to coincide with the $a_i^j$ and makes the two sections $s$ and $t$ coincide in the neighborhood of these points. 

The other zeroes of $s^1$ of index $i$ %intersection points of $s$ with $\mathcal M_0^i$  
can be matched  into pairs of points $\{a_i^{j+},a_i^{j-}\}$ of opposite sign. 
A Whitney type lemma allows us  to cancel all these pairs by a suitable homotopy 
of $s$ in $\mathcal M$, reducing to the case when $s^1$ has no other zeroes than 
the $a_i^j$'s, $j=1,\ldots,\vert m_i\vert$. A similar reduction may be assumed for $t$.
 Let us finish the proof in this case before stating and proving this lemma.

Both sections $s^1$ and $t^1$ of $T^*M$ are homotopic (among sections) by a homotopy which is stationary on a neighborhood $N(a_i^j)$. Making this homotopy 
$h: M\times[0,1]\to T^*M$ transverse to the zero section, the preimage of $0_M$ consists of arcs
$\{a_i^j\}\times[0,1]$ and finitely many closed curves $\ga_k$. Each of these closed curves can be arbitrarily decorated with an index $i$.  This choice allows us to lift
$h$ to $\widetilde J^2(M)$ %and 
as a  homotopy $\tilde h$ from $s$ to $t$; this 
$\tilde h$ %being 
is the desired homotopy. More precisely, we proceed as follows for getting 
$\tilde h$. First $h\vert \ga_k$ is lifted to $\mathcal M_0^i$ by using that the fiber of $\pi: \mathcal M_0^i\to 0_M$ is connected. The transversality of $h$ to $0_M$ allows us to extend this lifting to a neighborhood of $\ga_k$, making 
$\tilde h$ transverse  to $\mathcal M_0^i$. Now it is easy to extend  $\tilde h$ to $M\times [0,1]$, since 
the fiber of $\pi$ over any point outside $0_M$ is contractible.\\

\nd {\sc A Whitney type lemma.} {\it Let $(b^+,b^-)$ be a pair of transverse 
intersection points of $s$ with $\mathcal M^i_0$
having  opposite sign when they are thought of as zeroes of $s_1$ in $M$. 
%(seen in $M$, the source of $s$, as zeroes of $s^1$). 
Let $\al$ be a simple path in $M$ joining them 
avoiding the other zeroes of $s^1$ and let $N$ be a neighborhood of $\al$. Then there exists
a homotopy $S=(s_u)_{u\in[0,1]}$ of $s_0=s$ into $\mathcal M$, supported in $N$ 
and cancelling the pair $(b^+,b^-)$, that is, $\pi\circ s_1$ has no zeroes in $N$.}\\

\nd {\sc Proof.} We choose an embedded
 2-disk (with corners) $\De$ in $N\times [0,1[$ meeting $N\times \{0\}$ transversally 
 along $\al $. We first construct the homotopy $S^1:=\pi\circ S$ of 
 $s^1$ among the sections of $T^*M$, following the cancellation
 process of Whitney which we are going to recall. We require $S^1$ to be transverse to $0_M$ with 
 $(S^1)^{-1}(0_M)=\beta$, where $\alpha\cup\beta= \partial\De$. 
 Using a trivialization of $T^*M\vert N$, $S^1\vert N\times [0,1] $
  reads $S^1(x,u)=(x,g(x,u))$. The requirement is that $g$ vanishes transversally along
  the arc $\beta$; it is possible exactly because $\dim M\geq 2$ and 
  the end points have opposite signs.
  Let $T$ be a small tubular neighborhood of $\beta$; its boundary traces an arc 
  $\beta'$ on $\De$, ``parallel'' to $\beta$. Let $\al'$ be the subarc of $\al$ whose end points are those of $\beta'$. The restriction $g\vert T$ is required to be a trivialization
  of $T$, but this latter may be chosen freely. We choose it so that the loop
  $(g\vert\beta')\cup(s^1\vert\alpha')$ be homotopic to 0 in $(\R^n)^*\setminus\{0\}$; 
  of course, when $n>2$ this condition is automatically fulfilled. Now $g$ can be extended
  to the rest of $\De$ as a non-vanishing map. As $N\times[0,1]$ collapses onto
  $N\times\{0\}\cup \De\cup T$, the extension of $g$ can be completed without adding
 zeroes outside $\beta$, yielding the desired homotopy $S^1$.
 
 It remains to lift $S^1$ to $\mathcal M$. The lifting is first performed along $\beta$
 with value in $\mathcal M_0^i$. Then it is globally extended in the same way
 as in the above lifting process.
\bull\\

\begin{cor} Let $s$ be a section of $\mathcal M\subset\widetilde J^2(M)$ and $m_i$ be its algebraic intersection number with $\mathcal M_0^i$. 
Let $f:M\to\R$ be a Morse function whose number $c_i(f)$ of critical points of index $i$
satisfies  $$ c_i(f)= m_i$$ for all $i\in \{0, \ldots, n\}$. Then $s$ and $j^2(f)$ are homotopic as sections of 
 $\mathcal M$.\label{homot-cor}
\end{cor}

\begin{cor} We assume $\dim M\geq 6$ and $\pi_1(M)=0$. Let $s$ be a section of 
$\mathcal M\subset\widetilde J^2(M)$ whose algebraic intersection numbers
$m_i$ fulfills the Morse inequalities for every field of coefficients.
In particular, they are non-negative. Then $s$ is homotopic through sections in 
$\mathcal M$ to a holonomic section.
\end{cor}

\proof Under these topological assumptions  the following result holds true:
{\it For any set of non-negative integers $\{c_0,c_1,\cdots, c_n\}$ satisfying the Morse 
inequalities
for any field of coefficients,  there exists a Morse function on $M$ with $c_i$ critical
 points of index $i$} (see theorem 2.3 in \cite{franks}). So we have a Morse function 
$f:M\to \R$ with $m_i$ critical points of index $i$. According to corollary
\ref{homot-cor}, $s$ is homotopic in $\mathcal M$ to $j^2f$.\bull\\

\begin{rien}{\rm We end this section by recalling that the Morse inequalities are not sharp for estimating the number of critical points of a Morse function 
on a non-simply connected closed manifold. 
Typically when $\pi_1(M)$ equals its subgroup of commutators (perfect group), some critical points of index 1 are required for generating the fundamental group, but the Morse inequalities allow $c_1=0$ (see \cite{maller} 
for more details). On the other hand, the only constraint for a section of  $\mathcal M$ with intersection numbers $m_i$ is the Euler-Poincar\'e identity:
$$m_0-m_1+...=\chi(M).
$$ So it is possible to find a section $s$ whose intersection number $m_i$ is  the minimal rank in degree $i$ of a free complex whose homology is $H_*(M,\Z)$, that is,
$$m_i=\beta_i+\tau_i+\tau_{i-1}, $$ where $\beta_i$ stands for the rank of the free
quotient of $H_i(M,\Z)$ and $\tau_i$ denotes the minimal number of generators of its torsion subgroup (\cite{franks} p. 15). Such a set of integers satisfies the Morse inequalities but is far from being realizable by a Morse function. Finally this section $s$ is not homotopic in 
$\mathcal M$ to a holonomic section.}\\
\end{rien}

\section{Failure of the 1-parametric version of the $h$-principle}

\medskip
We thank Yasha Eliashberg who pointed out to us the failure of the $h$-principle in the
1-parametric version of the problem under consideration.

Here $M$ is assumed to be a product  $M=N\times [0,1]$. 
Let $f_0: M\to [0,1]$ be the projection. 
When $M$ is not 1-connected and $\dim M\geq 6$, according to Allen Hatcher 
the so-called {\it pseudo-isotopy 
problem} has always a negative answer: there exists $f$
without critical points which is not joinable to $f_0$ among the Morse functions (see \cite{hatcher}). But $j^2f$ can be joined to $j^2f_0$ by a path $\ga$ in 
$\mathcal M$. Indeed, take a generic homotopy $\ga^1$ joining $df$ to $df_0$; then arguing as in the proof of proposition \ref{homot} it is possible to lift it to 
$\mathcal M$. When $M$ is the $n$-torus $\T^n$, A. Douady showed very simply the stronger fact  that the path  $\ga^1$ can be taken among the non-singular 1-forms (see appendix to \cite{laud}).
 This $\ga$ is not homotopic in $\mathcal M$ with end points fixed to a path of holonomic sections.


\vskip 1cm
\begin{thebibliography}{99}
\bibitem{el-mi} Y. Eliashberg, N. Mishachev, {\it  Introduction to the  $h$-principle}, GSM 48, Amer. Math. Soc.,  2002.
\bibitem{franks} J. Franks, {\it Homology and Dynamical Systems}, CBMS Regional
Conf. vol. 49, Amer. Math. Soc., 1982.
\bibitem{gromov} M. Gromov, {\it Partial Differential Relations,} Springer Verlag, 1986.
 \bibitem{hatcher} A. Hatcher,   {\it Higher simple homotopy theory,} 
 Annals of Math. (2) 102 (1975), 101-137.
 \bibitem{laud}  F. Laudenbach, {\it Formes diff\'erentielles de degr\'e 1 ferm\'ees non singuli\`eres : classes d'homotopie de leurs noyaux}, Commentarii Math. Helvetici 
 51 n$^o$3 (1976), 447-464.
  \bibitem{bordeaux}  F. Laudenbach, {\it De la tranversalit\'e de Thom au $h$-principe de Gromov}, Le\c{c}ons de Math\'ematiques d'Aujourd'hui, vol. 4, \'Ed. Cassini, Paris (to appear).
  \bibitem{maller} M. Maller, {\it Fitted diffeomorphisms of non-simply connected manifolds, }Topology 19 (1980), 395-410.
  \bibitem{smale} S. Smale, {\it Notes on differentiable dynamical systems}, 277-287 
  in: Proc. Symposia Pure Math., vol. 14, Amer. Math. Soc., 1970.
  \bibitem{vassiliev} V.A. Vassiliev,{\it Topology of spaces of functions without compound singularities}, 23 (4) (1989), 277-286.
  
  \end{thebibliography}
\end{document}